\documentclass[12pt]{iopart}
\usepackage{iopams}
\newcommand{\nc}{\newcommand*}
\newcommand{\rnc}{\renewcommand*}
\rnc{\Bbb}{\mathbb}
\nc{\Cal}{\mathcal}

\def\BZ{\Bbb Z}
\def\CF{\Cal F}
\def\CR{\Cal R}
\def\CU{\Cal U}
\def\la{\lambda}
\def\vp{\varepsilon}

\def\Ker{\operatorname{Ker}}
\def\Ker{{\rm Ker}}

\begin{document}

\title[Deformation of $osp(1|4)$]
{Deformation of orthosymplectic Lie superalgebra $osp(1|4)$}

\author{Enrico Celeghini\dag and Petr P Kulish\ddag
\footnote[3]{Send the correspondence to P P Kulish (kulish@pdmi.ras.ru)}
}

\address{\dag\ Dipartimento di Fisica, Universita di Firenze
and INFN Sezione di Firenze,
Via G. Sansone 1, Sesto, 50019, Firenze, Italy}

\address{\ddag\
St. Petersburg Department of Steklov Mathematical Institute,
Fontanka 27, St. Petersburg 191023, Russia}

\begin{abstract}
Triangular deformation of the
orthosymplectic Lie superalgebra $osp (1|4)$ is defined by chains of
twists. Corresponding classical $r$-matrix is obtained by a contraction
procedure from the trigonometric $r$-matrix. The carrier space of the
constant $r$-matrix is the Borel subalgebra.
\end{abstract}


\submitto{\JPA}


%
%
%
%
%
%
%

\section{Introduction}

Orthosymplectic Lie superalgebras $osp(m|2n)$  have variety of applications
in gauged supergravity,  supersymmetric quantum mechanics, integrable
$\BZ_2$-graded spin chains etc. Similar applications with corresponding
changes can have quantum deformations of these superalgebras. In this paper
a  triangular deformation of the basic Lie superalgebra $g=osp(1|4)$ is
given by explicite construction of a universal twist element
$\CF \in \CU(g) \otimes \CU(g)$ \cite{1}.  
It is a natural extension of  known triangular deformation of rank 1 Lie
superalgebra $osp(1|2)$ \cite{2, 3}    
along the lines presented in \cite{4}    
about twists with deformed carrier subspaces. The constructed twist includes
all the generators  of the Borel subalgebra
$B_+\subset osp(1|4)$. Deformed
quantum algebra can be used to define a noncommutative version of the super
anti-de Sitter  space ($s$-AdS) generalizing \cite{5}--\cite{8},    
deformation of superconformal  quantum mechanics \cite{9},
deformed spin chains \cite{PK} etc.

\section{Classical $r$-matrices }

The direction of twist deformation is given by a classical  triangular
$r$-matrix $r\in g\wedge g$, which is a solution
of classical (super) Yang--Baxter
equation (cYBE) on $g\otimes g\otimes g$  \cite{10}
\begin{equation}\label{tag1}
[r_{12}, r_{13}+r_{23}] + [r_{13}, r_{23}]=0,
\end{equation}
where in the super-case the commutators and the tensor products are
$\BZ_2$-graded:
\begin{eqnarray}
[a,b]=ab -(-1)^{p(a)p(b)}ba, \nonumber \\
(a\otimes b)(c\otimes d)=(-1)^{p(b)p(c)} (ac\otimes bd).\nonumber
\end{eqnarray}
The parity of a homogeneous element $b$ of a $\Bbb Z_2$-graded space
 $V=V_0\oplus V_1$
is denoted by $p(b)$. For  basic classical Lie superalgebras
the classical $r$-matrix we are looking for, can be obtained by a contraction
procedure from the trigonometric $r$-matrix
\begin{equation}
r(\la -\mu)=(r_0e^{\la-\mu}+r^{21}_0)/ (e^{\la-\mu}-1),
\end{equation}
where $r_0$ is the super-analog of the Drinfeld--Jimbo constant $r$-matrix,
the solution to the cYBE,
\begin{equation}\label{r0}
r_0= \frac 12 \sum k_i \otimes k_i
+\sum\limits_{\alpha \in \Delta_+} e_{\alpha}\otimes e_{-\alpha}.
\end{equation}
Here $\{k_i\}$ is an orthonormal basis of the Cartan subalgebra and $\Delta_+$
is the set of positive (even and odd) roots. The sum
\[ 
c^{\otimes}_{2}:=r_0+r^{21}_0=\sum k_i\otimes k_i+
\sum\limits_{\alpha \in \Delta_+} (e_{\alpha}\otimes e_{-\alpha}+
(-1)^{p(\alpha)}e_{-\alpha}\otimes e_{\alpha})
\]
is an element of $\CU(g)\otimes \CU (g)$ (the tensor Casimir element)
invariant with respect to  the adjoint action:
\[ 
[x\otimes 1+1\otimes x,c^{\otimes}_2]=0, \quad  x\in g.
\]

Take the long root generator $e_{\theta}$
 and consider adjoint transformation
by the group element $\exp (te_{\theta})$:
\begin{eqnarray}
\fl  Ad (\exp(te_{\theta}))^{\otimes 2}r(\la -\mu)=
\frac 12\bigg(\coth \frac{\la-\mu}{2} \cdot c^{\otimes}_2+ \nonumber\\
 \sum_{\alpha \in \Delta_+}(e_{\alpha}\otimes e_{-\alpha}-
(-1)^{p(\alpha)}e_{-\alpha} \otimes e_{\alpha})+t\sum_{\alpha\in \Delta_+}
e_{\alpha}\wedge [e_{\theta},e_{-\alpha}]\bigg).
\end{eqnarray}
Scaling the spectral parameter $\la \to \vp \la$, $\mu \to \vp \mu$ and
$t\to 2t/\vp$, one gets after contraction
\begin{eqnarray}
\lim\limits_{\vp\to 0}
\vp Ad(\exp(t e_{\theta}/\vp))^{\otimes2}r(\vp(\la -\mu))=
\frac {c^{\otimes}_2}{(\la -\mu)} +
t\sum_{\alpha\in \Delta_+}e_{\alpha}\wedge
[e_{\theta},e_{-\alpha}].
\label{tag2}
\end{eqnarray}
Due to the Lie algebra nature of the cYBE (1),
the resulting expression satisfies
this equation as well and both terms satisfy it separately. The second term
will be used to deform the universal enveloping algebra $\CU (osp(1|4))$.
The Cartan--Weyl basis of the orthosymplectic Lie superalgebra $osp(1|4)$ has
the following generators:  $H,J$ (the Cartan subalgebra), $v_{\pm}$, $w_{\pm}$
(the odd root generators),
$X_{\pm}=\pm (v_{\pm})^2$, $Y_{\pm}=\pm(w_{\pm})^2$,
$U_{\pm}$, $Z_{\pm}$. There are two $osp(1|2)$ subalgebras $\{ H,v_{\pm},
X_{\pm}\}$ and $\{ J,w_{\pm}, Y_{\pm}\}$,  and the even
root generators define the Lie subalgebra $sp(4)\simeq so(5)$.
Constructing a universal twist element
$\CF$ we shall use only commutators of the Borel subalgebra generators with
subscripts $+$. Some of these commutation  relations are
\begin{equation}
[H,X_{+}] = 2X_+, \quad [H,v_+]=v_+,\quad [H,U_+]=U_+\,, 
\end{equation}
\begin{equation}
[H,Z_{+}] = Z_+, \quad [Z_+,U_+]=2X_+,\quad [Z_+,Y_+]=U_+\,, 
\end{equation}
\begin{equation}
[Z_+,w_{+}] = v_+, \quad [v_+,w_+]=U_+\,. 
\end{equation}
The generators $Z_+$, $w_+$ correspond to the simple roots. $X_+$ corresponds
to the long root $\theta$ and
it commutes with all positive root generators and with
the generator $J$ of the Cartan subalgebra.  The root $\theta$ in (4), (5) 
is the long one, hence all generators in the constant $r$-matrix (5) 
correspond to positive roots and one generator belongs to the Cartan subalgebra
$[X_{\theta},X_{-\theta}]=H$.  In the case of $osp(1|4)$ this constant
$r$-matrix is $(X_{\theta}=X_+)$
\begin{equation}\label{esjr}
r=H\wedge X_+-v_+\otimes v_++Z_+\wedge U_+.
\end{equation}
This is an extension by the odd generator $v_+$ of the classical $r$-matrix
defining the extended jordanian twisting of the
$\CU (so(5)) \simeq \CU (sp(4))$ \cite{4,17}.

Solution of the cYBE defines a cobracket $\delta$ on the Lie (super) algebra:
a map $\delta : g \to g \wedge g$,
\[
\delta (x)=[x\otimes 1+1\otimes x,r],\quad x\in g.
\]
It is important to point out that the dimension of the kernel of $\delta$ is
equal to the number of independent primitive elements in the quantum
$\CU_{\xi}(g)$ which is a quantization of the bialgebra $\CU(g)$. Proof
can be given by transition to the dual Poisson--Lie group using
the quantum duality principle \cite{11, 12}.   

The kernel of cobracket defined by the $r$-matrix (\ref{esjr})   
is generated by the elements $X_+$, $J$, $Y_+$, $w_+$.
It is easy to see that
the kernel of the cobracket $\delta$ is a Lie subalgebra $\Ker \delta
\subset  g$. If there is a triangular $r$-matrix with its carrier in
$\Ker \delta$, then the sum of this additional $r$-matrix and the one defined
$\delta$, satisfies the cYBE as well. In the case we are interested in there
is a super-jordanian $r$-matrix \cite{3}  
\begin{equation}\label{tag4}
r^{(sj)}=J\wedge Y_+-w_+\wedge w_+.
\end{equation}
Finally, we shall construct a universal twist corresponding to the $r$-matrix:
\begin{equation}\label{esj2}
r=H\wedge X_+-v_+\otimes v_+ +Z_+\wedge U_+ +J\wedge Y_+-w_+\otimes w_+.
\end{equation}
(One can add also an abelian term $X_+\wedge Y_+$, and introduce few
independent parameters in front of different constituent $r$-matrices.)

This form of the  classical $r$-matrix can be generalized
to the case of the universal enveloping algebra of the orthosymplectic Lie
superalgebra $osp(1|2n)$. It has  (positive) even simple root generators
$Z_k$, $k=1,2,\dots,n-1$, $\alpha_k=\varepsilon_k-\varepsilon_{k+1}$,  and
one odd simple root generator $v_n$, $\alpha_n=\varepsilon_n$. Other positive
root generators are $X_j=v^2_j$, $j=1,2,\dots,n$, $Z_{kj}$ and $U_{kj}$ related
to positive roots $\varepsilon_k-\varepsilon_j$ and $\varepsilon_k+
\varepsilon_j$, $1\le k<j\le n$ \cite{13}. Commutation relations
\begin{equation}\label{12}
[H_k,X_k]=2X_k, \quad [H_k,v_k]=v_k,
\end{equation}
\begin{equation}\label{13}
[Z_{kj}, U_{kj}] =2X_k, \quad [v_k,v_j]=U_{kj},
\end{equation}
\begin{equation}\label{14}
[Z_{kj},v_j]=v_k,\quad [H_k,Z_{kj}]=Z_{kj}, \quad [H_k,U_{kj}]=U_{kj}
\end{equation}
are used to prove that
\begin{equation}\label{2nr}
r=\sum^{n}_{k=1}a_k(H_k\wedge X_k+\sum^{n}_{j>k}Z_{kj}\wedge U_{kj}
-v_k\otimes v_k)
\end{equation}
satisfies the cYBE. Similar to the symplectic algebra $sp(2n)$ \cite{18}
the sum of $r$-matrices (\ref{2nr}) corresponds to the sequence of injections
\[
osp(1|2)\subset osp(1|4)\subset \dots \subset osp (1|2n).
\]

\section{Universal twist element of $B_+\subset osp(1|4)$ }

Construction of an explicit form of the twist corresponding to the
triangular $r$-matrix  (\ref{esjr})
will be realized according to the recipe
known as chains of twists \cite{14, 15}.  
The classical $r$-matrix is decomposed into few terms which generate factors
of the universal twist. This decomposition consists of the jordanian $r$-matrix
\begin{equation}\label{tag7}
r^{(j)}=H\otimes X_+-X_+\otimes H:=H\wedge X_+,
\end{equation}
extended jordanian $r^{(ej)}$ and super-jordanian $r$-matrices
\begin{equation}\label{tag8}
r^{(ej)}=r^{(j)}+Z_+\wedge U_+,\quad r^{(sj)}=r^{(j)}-v_+\otimes v_+.
\end{equation}
Twist elements corresponding to $r^{(ej)}$ and $r^{(sj)}$  are
known \cite{14, 3}, and they are represented in a factorized form
\begin{equation}\label{tag9}
\CF^{(ej)}=\CF^{(e)}\CF^{(j)}, \quad \CF^{(sj)}=\CF^{(s)}\CF^{(j)},
\end{equation}
where
\begin{equation}\label{19}
\CF^{(j)}=\exp (H\otimes \sigma), \quad \CF^{(e)}=\exp (\frac 12Z_+
\otimes U_+e^{-\sigma}),
\end{equation}
\begin{equation}\label{20}
\CF^{(s)}=\exp (-v_+\otimes v_+f(\sigma \otimes 1,1\otimes \sigma)),\quad
\sigma =\frac 12\ln (1+X_+).
\end{equation}
The function $f$ is symmetric. It can be presented as the series expansion
$\sum f_n$ \cite{3}, or through the factorization of $\CF^{(s)}$ as
\begin{equation}\label{21}
\CF^{(s)}=
(1-(v_+\otimes v_+)(f_1(\sigma) \otimes f_1(\sigma)))\CF^{(c)}\;,
\qquad f_1(\sigma)=(e^{\sigma}+1)^{-1}
\end{equation}
with an appropriate coboundary twist \cite{9}
\[
\CF^{(c)}=(u\otimes u)\Delta (u^{-1})\;,\qquad
u=(\frac12(e^{\sigma}+1))^{\frac12}\;.
\]
Due to the commutativity of $X_+$ or $\sigma$ and
$v_+$ with $Z_+$, $U_+$, the
twisting by $\CF^{(s)}$ after $\CF^{(j)}$ does not change the coproducts
of $Z_+$ and $U_+$, and vise versa.  Hence, one can arrange these twists
to form an extended superjordanian  twist
\begin{equation}\label{esjF}
\CF^{(esj)}=\CF^{(s)}\CF^{(e)}\CF^{(j)}=\CF^{(e)}\CF^{(s)}\CF^{(j)},
\end{equation}
corresponding to the $r$-matrix (\ref{esjr}).

To take into account further extension of the $r$-matrix (\ref{esj2})
by a super-jordanian term of the second $osp(1|2)$ subalgebra, we have  to use
a deformed carrier space transforming the generators $Y_+$ and $w_+$ according
to \cite{4, 14}.

The generators $v_+, U_+$ and $\sigma (X_+) = \frac12 \ln(1+X_+)$
are entering into this transformation. After the twist (\ref{esjF})
their coproducts are
\[
\Delta^{(esj)}(\sigma) := \CF^{(esj)} \Delta(\sigma) (\CF^{(esj)})^{-1} =
\Delta_0(\sigma) = \sigma \otimes 1 + 1 \otimes  \sigma \;,
\]
\[
\Delta^{(esj)}(v_+) = \Delta^{(sj)}(v_+) = v_+ \otimes 1 + e^{\sigma} \otimes
v_+ \;,
\]
\[
\Delta^{(esj)}(U_+) = \Delta^{(ej)}(U_+) = U_+ \otimes e^{\sigma} +
e^{2\sigma} \otimes U_+ \;.
\]

The twisted coproducts of the Borel subalgebra generators $J,Y_+,w_+$ of the
second  $osp(1|2)$ are
\[
\Delta^{(esj)}(J) = \CF^{(esj)} \Delta(J) (\CF^{(esj)})^{-1} 
= \Delta_0(J) = J \otimes 1 + 1 \otimes  J \;,
\]
\[
\Delta^{(esj)}(Y_+) = \Delta_0(Y_+) +\frac12 U_+ \otimes U_+  e^{-\sigma} +
\frac14 X_+ \otimes U_+^2  e^{-2\sigma}\;,
\]
and we get a rather lengthly expression for $\Delta^{(esj)}(w_+)$.

The deformed generator $\widetilde Y_+$ with the primitive coproduct is
given by the adjoint transformation \cite{4,18}
\begin{equation}\label{Ad}
\widetilde Y_+= Ad \exp (-Z_+U_+\sigma /2X_+)\,Y_+ =
Y_+-\frac 14U^2_+e^{-2\sigma}.
\end{equation}
Hence, the second jordanian factor is similar to the $sp(4)\simeq so(5)$ case
\cite{4}  and has the form $\exp (J\otimes \sigma (\widetilde Y_+))$.
Although, the deformed
coproduct of the odd generator $\widetilde{w}_+$ is rather cumbersome after the
transformation (\ref{Ad}) the coproduct of $\widetilde{w}_+$
\[
\widetilde{w}_+ =  Ad \exp(-Z_+U_+\sigma/2X_+) \, w_+=w_+-
\frac12 v_+U_+\frac{e^{-\sigma}}{e^{\sigma}+1}
\]
is primitive as well. Hence, one can add to (\ref{esjF}) an additional
factor $\CF^{(sj2)}$ corresponding to the second  $osp(1|2)$
\[
\CF^{(sj2)} = \exp (-\widetilde{w}_+ \otimes \widetilde{w}_+
f(\widetilde\sigma \otimes 1,1\otimes \widetilde\sigma))
\exp (J\otimes \sigma (\widetilde Y_+)) \,.
\]
where $\widetilde\sigma :=\sigma (\widetilde Y_+) 
= \frac 12 \ln (1 + \widetilde Y_+) $. 

A universal twist with the carrier space including all generators
of the Borel subalgebra $B_+\subset osp(1|4)$, is  given by the product
\[
\CF=\CF^{(sj2)}\CF^{(esj)}\;,
\]
where the super-jordanian twist $\CF^{(sj2)}$ \cite{3} is constructed using
the generators $J, \widetilde Y_+, \widetilde{w}_+$ as in (18), (20).
The universal $R$-matrix \cite{1} of the twisted $\CU(osp(1|4))$ is
\[
\CR = \CF_{21} \, \CF^{-1}\,.
\]

It is natural to conjecture that a similar chain of twists  can be
constructed for $\CU(osp(1|2n))$ starting with the $r$-matrix (\ref{2nr})
as it was in the case of $\CU(sp(2n))$ and
$\CU(so(2n+1))$ \cite{18, 19}.  
This would extend in some sense  an analogy between algebras $so(2n+1)$
and $osp(1|2n)$ \cite{16}.     

One can introduce more convenient set of generators of the twisted Borel
subalgebra $B_+\subset osp (1|4)$ using the $FRT$-formalism \cite{20}, and
the upper triangular form of the generators of $B_+$ in the defining
$5\times 5$ irreducible representation $\rho$. The new generators are
the entries of
the universal $L$-operator: $L=(\rho \otimes id)\Cal R$.
Using the symmetric
grading $(0,0,1,0,0)$ of the rows and  columns one gets
\[
L=\left(
\begin{array}{ccccc}
   T^{-1}&u     & V&z          &h\\
        0&K^{-1}&W &j           &\tilde z \\
        0&0     &1 &\widetilde W&\widetilde V\\
        0&0     &0 &K           &\tilde u\\
        0&0     &0 &0           &T
\end{array}
\right)\,.
\]
The coproduct of these generators is given by the matrix product of the
$L$-operators. The commutation relations follow from the $RTT$-relation:
$RL_1L_2=L_2L_1R$ \cite{20},  taking into account extra signs in the
$Z_2$-graded tensor
product \cite{10}. The $R$-matrix can be obtained using
also a super-version of
the $r^3=0$ theorem: $R=\exp (\eta r_{\rho})$, where
$r_{\rho}=(\rho\otimes \rho)r$. Similarly to the
twisting of $osp(1|2)$ \cite{3} one
can express the new generators in terms of the initial ones.

\subsection*{Acknowledgement}
We are grateful to V N Tolstoy for useful discussion and
the information that a similar result was
obtained in \cite{21}.  We thank V D Lyakhovsky for valuable comments.
One of the authors (PPK) would like to thank the
INFN Sezione di Firenze for generous hospitatlity.
This work has been partially supported by the grant  RFBR-02-01-00085 and
the programme ``Mathematical methods in nonlinear dynamics'' of RAN.

\bigskip
{\bf References}
\bigskip

\def\cmp#1#2#3{ {#2} {\it Commun. Math. Phys.} {#1} {#3}}
\def\jpa#1#2#3{ {#2} {\it  J. Phys. A: Math. Gen.} {#1} {#3}}
\def\jmp#1#2#3{ {#2} {\it  J. Math. Phys.} {#1} {#3}}
\def\plb#1#2#3{ {#2} {\it  Phys. Lett. B.} {#1} {#3}}
\def\prd#1#2#3{ {#2} {\it  Phys. Rev. D.} {#1} {#3}}
\def\mpla#1#2#3{ {#2} {\it Mod. Phys. Lett. A.} {#1} {#3}}

\end{document}